\documentclass[12pt,english]{amsart}
\usepackage[T1]{fontenc}
\usepackage[latin9]{inputenc}
\usepackage{geometry}
\usepackage{amsthm}
\usepackage{amstext}
\usepackage{amssymb}
\usepackage{amsmath}
\usepackage{setspace}
\usepackage{esint}
\usepackage{enumerate}
\usepackage{etoolbox}
\usepackage{url}
\usepackage{mathrsfs}
\usepackage{graphicx}
\usepackage[labelfont = bf]{caption}  
\usepackage{amsaddr}
\usepackage{hyperref}
\hypersetup{
  colorlinks = false,
  pdfborder = {0,0,0}
}

\makeatletter
\numberwithin{equation}{section}
\numberwithin{figure}{section}
\theoremstyle{plain}
\newtheorem{thm}{\protect\theoremname}[section]
  \theoremstyle{plain}
  \newtheorem{prop}[thm]{\protect\propositionname}
  \theoremstyle{plain}
  
  \theoremstyle{definition}
  \newtheorem{defn}[thm]{\protect\definitionname}
  \theoremstyle{remark}
  \newtheorem{rem}[thm]{\protect\remarkname}
  \theoremstyle{plain}
  \newtheorem{lem}[thm]{\protect\lemmaname}
  \theoremstyle{definition}

\makeatother

\usepackage{babel}
  \providecommand{\corollaryname}{Corollary}
  \providecommand{\definitionname}{Definition}
  \providecommand{\lemmaname}{Lemma}
  \providecommand{\propositionname}{Proposition}
  \providecommand{\remarkname}{Remark}
  \providecommand{\examplename}{Example}
\providecommand{\theoremname}{Theorem}

\patchcmd{\section}{\scshape}{\bfseries\large}{}{}
\makeatletter 
\renewcommand{\@secnumfont}{\bfseries\large}  
\makeatother

\makeatletter
{}{}
\makeatother

\makeatletter
\patchcmd{\@maketitle}
  {\ifx\@empty\@dedicatory}
  {\ifx\@empty\@date \else {\vskip3ex \centering\footnotesize\@date\par\vskip1ex}\fi
   \ifx\@empty\@dedicatory}
  {}{}
\patchcmd{\@adminfootnotes}
  {\ifx\@empty\@date\else \@footnotetext{\@setdate}\fi}
  {}{}{}
\makeatother

\renewcommand{\mathbb}{\mathbf}

\begin{document}


\title[Sharp boundary regularity for SPDEs in a half-line]{Sharp regularity near an absorbing boundary for solutions to second
order SPDEs in a half-line with constant coefficients}

\author{Sean Ledger}

\begin{abstract}  
We prove that the weak version of the SPDE problem
\begin{align*}
dV_{t}(x) & = [-\mu V_{t}'(x) + \frac{1}{2} (\sigma_{M}^{2} + \sigma_{I}^{2})V_{t}''(x)]dt - \sigma_{M} V_{t}'(x)dW^{M}_{t}, \quad x > 0,
\\
V_{t}(0) &= 0
\end{align*}
with a specified bounded initial density, $V_{0}$, and $W$ a standard
Brownian motion, has a unique solution in the class of finite-measure
valued processes. The solution has a smooth density process which
has a probabilistic representation and shows degeneracy near the absorbing
boundary. In the language of weighted Sobolev spaces, we describe
the precise order of integrability of the density and its derivatives
near the origin, and we relate this behaviour to a two-dimensional
Brownian motion in a wedge whose angle is a function of the ratio
$\sigma_{M}/\sigma_{I}$. Our results are sharp: we demonstrate that
better regularity is unattainable. 
\end{abstract}

\address{Mathematical Institute, University of Oxford, Woodstock Road, Oxford, OX2 6GG}
\email{ledger@maths.ox.ac.uk}

\maketitle


{\setstretch{1}
\tableofcontents{}
}

\section{Introduction} 

Let $W$ be a standard Brownian motion on a complete filtered probability
space $\left(\Omega,\mathcal{F},\mathbb{P}\right)$ and $V_{0}$ be
a deterministic bounded probability density function on the half-line
$\left(0,\infty\right)$. In this paper, we consider the problem of
finding a finite-measure valued process, $\nu=\left(\nu_{t}\right)_{t\in\left[0,T\right]}$,
satisfying the weak stochastic partial differential equation 
\begin{equation}
\left\langle \phi,\nu_{t}\right\rangle =\left\langle \phi,V_{0}\right\rangle +\int_{0}^{t}[\mu\left\langle \phi',\nu_{s}\right\rangle +\frac{1}{2}(\sigma_{M}^{2}+\sigma_{I}^{2})\left\langle \phi'',\nu_{s}\right\rangle ]ds+\sigma_{M}\int_{0}^{t}\left\langle \phi',\nu_{s}\right\rangle dW_{s}\label{eq:TheWeakProblem}
\end{equation}
for all $t\in\left[0,T\right]$ and $\phi\in C^{\textrm{test}}$.
Here, we use the abbreviation $\left\langle \phi,\nu_{t}\right\rangle =\int_{0}^{\infty}\phi(x)\nu_{t}(dx)$
and write $C^{\textrm{test}}$ for the space of bounded, twice-differentiable
functions that vanish at the origin and have bounded first and second
derivatives. The coefficients $\mu\in\mathbb{R}$, $\sigma_{M}>0$,
and $\sigma_{I}>0$ and the time horizon $T>0$ are deterministic
constants. We prove this problem has a unique solution, show how to
construct this solution, and establish its regularity near the absorbing
boundary.

This stochastic partial differential equation (SPDE) arises naturally as the mean-field limit of a collection of interacting particles.
Specifically, suppose we have a homogeneous pool of particles each
experiencing an independent noise driven by a Brownian motion with
speed $\sigma_{I}$ and drift $\mu$, as well as a common noise $\sigma_{M}W_{t}$.
If we kill each of the particles upon hitting the origin, then $\nu_{t}$
describes the limiting spatial distribution of the particles at time
$t$ as the number of particles becomes large. In \cite{Bush}, this
model is proposed for the pricing of portfolio credit derivatives,
and in that context the common noise, $W$, corresponds to a market
risk factor and the event that a particle hits the origin represents
a default of one of the names in the portfolio. As the particles in
this model are independent and identically distributed conditional
on knowing the history of $W$, we can study the system by considering
the conditional mean behaviour of a single particle. In Section 2
we use this idea to present a simple construction of the solution
to (\ref{eq:TheWeakProblem}).

Visualising the dynamics of the solution to (\ref{eq:TheWeakProblem})
is straightforward, and if $\sigma_{M}=0$ then the SPDE becomes the
deterministic heat equation with zero boundary condition. For non-zero
$\sigma_{M}$ the solution still has this decay, but also experiences
random fluctuations in location which are driven by the Brownian motion
$W$ --- see Figure 1. If we denote the density of the solution by
$V$, then we can write (\ref{eq:TheWeakProblem}) in differential
form: 
\[
dV_{t}=[-\mu V_{t}'+\frac{1}{2} (\sigma_{M}^{2}+\sigma_{I}^{2}) V_{t}'']dt-\sigma_{M}V_{t}'dW_{t},
\]
and, thinking of this equation as a PDE with irregular coefficients,
we see that $V$ diffuses with speed $\sigma_{M}^{2}+\sigma_{I}^{2}$
and the transport term, $\mu dt+\sigma_{M}dW_{t}$, describes the
(stochastic) net movement of the density in space. Lemma 3.1 of \cite{KrylovHeat}
shows that the density of the solution is infinitely differentiable
in $\left(0,\infty\right)$, and this is due to the smoothing effect
that heat flow produces. However, the presence of the irregular fluctuations
in the density's location causes mass close to the boundary to be
destroyed in an abrupt fashion, and as a result the derivatives of
the density blow-up at the origin. The purpose of this paper is to
describe this degeneracy at the absorbing boundary. 

\begin{figure}
\begin{centering}
\includegraphics[width=14cm,height=6cm, bb=0bp 60bp 576bp 315bp,clip]{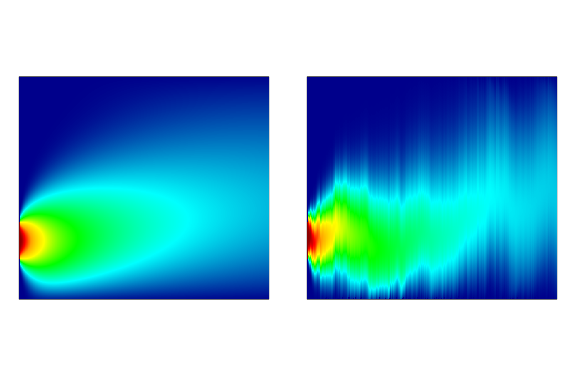}
\end{centering}
{\setstretch{1.0}
\caption{{\smaller  On the left is the evolution of the system in the deterministic case
$\sigma_{M}=0$. The value of the density at a space-time point is
represented by a colour increasing from blue ($V_{t}\left(x\right)=0$)
through to red ($V_{t}\left(x\right)=1$). The initial condition is
the step function $V_{0}=\mathbf{1}_{\left(0.18,0.5\right)}$ and
the drift, $\mu$, is zero. Compare this with a realisation of the
evolution for same parameters but with $\sigma_{M}>0$, which is presented
on the right. We still observe heat-type dispersion but with a non-smooth
fluctuation in the density's location.}}
}
\end{figure}

Dirichlet boundary problems for SPDEs of this type have been studied
extensively --- for example \cite{Kim, Krylov, KryLoto}.
As we are working in one-dimension with constant coefficients, the
simplicity of our setting allows greater insight into the behaviour
of the solution near the absorbing boundary. We take advantage of
a straightforward smoothing technique to present a more elementary
study of the problem than is currently offered in the standard literature,
and we are able to prove uniqueness of the solution to (\ref{eq:TheWeakProblem})
in the broad class of finite-measure valued processes. Throughout
this paper, when referring to \emph{uniqueness} we shall always mean
that if $\nu^{1}$ and $\nu^{2}$ are two finite-measure valued processes
solving (\ref{eq:TheWeakProblem}), then 
\[
\mathbb{E}\int_{0}^{T}\left(\left\langle \phi,\nu_{t}\right\rangle -\left\langle \phi,\bar{\nu}_{t}\right\rangle \right)^{2}dt=0,\quad\mathrm{for\, every}\,\phi\in C^{\textrm{test}}.
\]

The results closest to those in this paper are found in \cite{Krylov};
there, Krylov demonstrates the existence of a unique solution to (\ref{eq:TheWeakProblem})
in the class of twice-differentiable function-valued processes, $U=\left(U_{t}\left(x\right)\right)_{t\in\left[0,T\right]}$,
that satisfy 
\[
\mathbb{E}\int_{0}^{T}\bigl(\left\Vert U_{t}\right\Vert _{L^{2}\left(0,\infty\right)}^{2}+\left\Vert U_{t}'\right\Vert _{L^{2}\left(0,\infty\right)}^{2}+\left\Vert xU_{t}''\right\Vert _{L^{2}\left(0,\infty\right)}^{2}\bigr)dt<\infty.
\]
Furthermore, Krylov shows that if the initial density is $N$-times
differentiable with the function $x\mapsto x^{k}V_{0}^{\left(k\right)}\left(x\right)$
in $L^{2}\left(0,\infty\right)$ for $k=0,1,\ldots,N$, then the solution
is $\left(N+1\right)$-times differentiable with 
\begin{equation}
\mathbb{E}\int_{0}^{T}\bigl\Vert x^{k-1}U_{t}^{\left(k\right)}\bigr\Vert _{L^{2}\left(0,\infty\right)}^{2}dt<\infty\label{eq:KrylovsReg}
\end{equation}
for $k=0,1,\ldots,N+1$. (For a function $f$ on $\left(0,\infty\right)$,
we will always denote the $n^{\mathrm{th}}$ order derivative of $f$
by $f^{\left(n\right)}$.)

A key ingredient in \cite{Krylov}, and the majority of the literature
on this problem, is the introduction of weighted Sobolev spaces. These
spaces enable us to quantify the rate at which the solution and its
derivatives blow-up near the origin. In our one-dimensional setting,
it suffices to use the functions 
\begin{equation}
w_{c}\left(x\right):=x^{c}e^{-x},\quad x>0,\label{eq:WeightingDefinition}
\end{equation}
for $c\in\mathbb{R}$. For a given order of derivative of the solution
to (\ref{eq:TheWeakProblem}), we can ask: how small can $c$ be taken
in the above weighting function whilst preserving square-integrability
of the product of that derivative and $w_{c}$? From (\ref{eq:KrylovsReg}),
we can certainly take any $c\geq k-1$ for the $k^{\mathrm{th}}$
derivative of the solution. A complete answer to this question is
provided by the two proceeding theorems, and a critical quantity is
\begin{equation}
\alpha:=\frac{\pi}{2}+\arcsin\rho\in\left(\pi/2,\pi\right),\quad\textrm{where }\,\rho:=\frac{\sigma_{M}}{\sqrt{\sigma_{M}^{2}+\sigma_{I}^{2}}}\in\left(0,1\right).\label{eq:alpha}
\end{equation}
In the context of credit derivatives, $\rho$ is referred to as the
\emph{correlation} in the portfolio and describes the tendency of
defaults to occur simultaneously. The closer $\rho$ is to one, the
larger the value of $\alpha$, and therefore, observing that $\pi/\alpha-1\in\left(0,1\right)$,
the following result states that \emph{higher correlation produces
less regularity near the absorbing boundary.}

\begin{thm}[Uniqueness and regularity]
\label{Reg}If $V_{0}$ is bounded, then there exists a unique solution, $\nu$,
to (\ref{eq:TheWeakProblem}) in the class of finite-measure valued
processes, and, for almost every $\left(\omega,t\right)\in\Omega\times\left[0,T\right]$,
$\nu_{t}$ has a density $V_{t}$ on $\left(0,\infty\right)$.

Furthermore, suppose that $V_{0}$ is $N$-times weakly differentiable
in $\left(0,\infty\right)$, and that for $k=0,1,\ldots,N$ 
\[
\bigl\Vert w_{k-\beta/2}V_{0}^{\left(k\right)}\bigr\Vert _{L^{2}\left(0,\infty\right)}<\infty,\quad\textrm{for every }\beta\in\left(-\infty,\pi/\alpha-1\right),
\]
where $w_{\left(\cdot\right)}$ and $\alpha$ are defined in (\ref{eq:WeightingDefinition})
and (\ref{eq:alpha}). Then, for almost all $\left(\omega,t\right)\in\Omega\times\left[0,T\right]$,
$V_{t}$ is $\left(N+1\right)$-times weakly differentiable in $\left(0,\infty\right)$,
and for $k=0,1,\ldots,N+1$
\[
\mathbb{E}\int_{0}^{T}\bigl\Vert w_{k-1-\beta/2}V_{t}^{\left(k\right)}\bigr\Vert _{L^{2}\left(0,\infty\right)}^{2}dt<\infty,\quad\textrm{for every }\beta\in\left(-\infty,\pi/\alpha-1\right).
\]
\end{thm}

This result shows that, at a time later than zero, the solution has
differentiability one order higher than the initial density, and is
one multiple of $x$ more regular at the origin (recall that $w_{k-\beta/2}\left(x\right)=x^{k-\beta/2}e^{-x}$).
We should notice that, because the initial density is integrable and
bounded, the condition $w_{-\beta/2}V_{0}\in L^{2}\left(0,\infty\right)$
is always satisfied. The hypotheses placed on the initial density
allow its derivatives to blow-up at the origin and grow sub-exponentially
towards infinity. We recover the results of \cite{Krylov} by taking
$\beta=0$ in Theorem \ref{Reg}, and therefore we have expanded upon
the known regularity of the solution near the absorbing boundary. 

Our results are sharp as the following result shows that we cannot
take $\beta>\pi/\alpha-1$ in Theorem \ref{Reg}:

\begin{thm}[Converse result]
\label{Converse}If $V_{0}$ satisfies the hypotheses of Theorem \ref{Reg}, then for
$k=0,1,\ldots,N+1$ we have 
\[
\mathbb{E}\int_{0}^{T}\bigl\Vert w_{k-1-\beta/2}V_{t}^{\left(k\right)}\bigr\Vert_{L^{2}\left(0,\infty\right)}^{2}dt=\infty,\quad\textrm{for every }\beta\in\left(\pi/\alpha-1,+\infty\right).
\]
\end{thm}

Since the tails of the weighting functions all decay at the same rate,
the only way we can obtain the blow-up observed in Theorem \ref{Converse}
is if the solution is degenerate near the absorbing boundary. With
this in mind, it is clear that the first derivative of the solution
cannot be bounded in a neighbourhood of the origin, otherwise we obtain
the contradiction 
\[
\int_{0}^{1}x^{-\beta}dx=\infty
\]
for $\beta\in\left(\pi/\alpha-1,1\right)\subseteq\left(0,1\right)$.
In the deterministic setting ($\sigma_{M}=0$) the first derivative
does not have this feature, and this has been remarked upon by Krylov
in \cite{KrylovHeat}. If we shift the density by setting 
\[
f\left(t,x\right):=V_{t}\left(x+\mu t+\sigma_{M}W_{t}\right),
\]
then, by the It\=o--Wentzell formula, we have that $f$ solves the deterministic
heat equation (with speed $\sigma_{I}$) in the random region 
\[
\left\{ \left(t,x\right):t\in\left[0,T\right],x>-\mu t-\sigma_{M}W_{t}\right\} 
\]
with Dirichlet boundary condition. Theorem 5.1 of \cite{KrylovHeat}
shows that there exists $\lambda\in\left(0,1\right)$ such that, with
probability one, there is a dense subset of $\left[0,T\right]$ on
which 
\[
\lim_{x\downarrow0}\frac{f\left(t,x\right)}{x^{\lambda}}=\infty.
\]
On this set it is clear that the density can have no spatial derivative
at the absorbing boundary.

In the next section we present the construction of the solution to
(\ref{eq:TheWeakProblem}) and describe its properties. A proof of
Proposition \ref{TDCE2} is presented in Section 3, and this is the
main probabilistic result used to prove the above two theorems in
Sections 4 and 5.

\subsubsection*{Acknowledgements}

I would like to thank Ben Hambly for introducing me to this problem
and related models, and Philippe Charmoy and James Leahy for helpful
conversations.


\section{The probabilistic solution}

We shall construct the solution to (\ref{eq:TheWeakProblem}) in this
section. We have not yet proved that (\ref{eq:TheWeakProblem}) can
admit only one solution, therefore we refer to the measure-valued
process presented here as the \emph{probabilistic solution}. 

\subsection*{Construction}

Let us introduce two independent standard Brownian motions, $W^{1}$
and $W^{2}$, that are also independent of $W$. We define two test
processes, $X^{1}$ and $X^{2}$, by 
\[
X_{t}^{i}:=X_{0}^{i}+\mu t+\sigma_{M}W_{t}+\sigma_{I}W_{t}^{i},
\]
where $X_{0}^{1}$ and $X_{0}^{2}$ are independent random variables
with common law $V_{0}$, and are also independent of every other
random variable introduced so far. To kill the particles upon hitting
zero, we set 
\[
\tau^{i}:=\inf\left\{ s\in\left[0,T\right]:X_{s}^{i}\leq0\right\} .
\]

Conditional on knowing the trajectory of $W$, $X^{1}$ and $X^{2}$
are independent and identically distributed, therefore we make the
following definition:

\begin{defn}[Probabilistic solution]
\label{ProbSolution} Let $\mathcal{F}^{W}=\left\{ \mathcal{F}_{t}^{W}\right\} _{t\in\left[0,T\right]}$
be the natural filtration of $W$. The \emph{probabilistic solution},
$\nu$, is defined to be the measure-valued process 
\[
\nu_{t}\left(S\right):=\mathbb{P}\left(\left.X_{t}^{1}\in S;t<\tau^{1}\right|\mathcal{F}_{t}^{W}\right)
\]
for all $t\in\left[0,T\right]$ and $S\subseteq\left(0,\infty\right)$
Lebesgue measurable. \end{defn}

 Throughout this paper ``$S\subseteq\left(0,\infty\right)$
measurable'' shall always mean that $S$ is Lebesgue measurable\emph{.}
The following is a useful representation:

\begin{prop}
\label{HittingProb2} For all $t\in\left[0,T\right]$ and $S\subseteq\left(0,\infty\right)$
measurable we have 
\[
\mathbb{E}\left[\nu_{t}\left(S\right)^{2}\right]=\mathbb{E}\left[\mathbb{E}\left[\left.\mathbf{1}_{X_{t}^{1}\in S,\, X_{t}^{2}\in S};t<\min\left(\tau^{1},\tau^{2}\right)\right|\mathcal{F}_{t}^{W}\right]\right] \!.
\]
\end{prop}
\begin{proof}
Since $X^{1}$ and $X^{2}$ are equal in
distribution, we have $\nu_{t}\left(S\right)=\mathbb{E}\left[\left.\mathbf{1}_{X_{t}^{i}\in S};t<\tau^{i}\right|\mathcal{F}_{t}^{W}\right]$
for $i=1$ and $2$. Hence 
\begin{align*}
\mathbb{E}\left[\nu_{t}\left(S\right)^{2}\right] 
  &= \mathbb{E}\left[\mathbb{E}\left[\left.\mathbf{1}_{X_{t}^{1}\in S};t<\tau^{1}\right|\mathcal{F}_{t}^{W}\right]^{2}\right] \\
   &= \mathbb{E}\left[\mathbb{E}\left[\left.\mathbf{1}_{X_{t}^{1}\in S};t<\tau^{1}\right|\mathcal{F}_{t}^{W}\right]\mathbb{E}\left[\left.\mathbf{1}_{X_{t}^{2}\in S};t<\tau^{2}\right|\mathcal{F}_{t}^{W}\right]\right] \\
 & = \mathbb{E}\left[\mathbb{E}\left[\left.\mathbf{1}_{X_{t}^{1}\in S,\, X_{t}^{2}\in S};t<\min\left(\tau^{1},\tau^{2}\right)\right|\mathcal{F}_{t}^{W}\right]\right]
\end{align*}
as required, where the final line is due to the conditional independence
of $X^{1}$ and $X^{2}$. 
\end{proof}

It is a simple calculation with It\=o's formula to show that $\nu$
is a solution of (\ref{eq:TheWeakProblem}):

\begin{prop}
\label{W.E.E.}There exists a full subset of $\Omega$ on which 
\[
\left\langle \phi,\nu_{t}\right\rangle =\left\langle \phi,V_{0}\right\rangle +\int_{0}^{t}[\mu\left\langle \phi',\nu_{s}\right\rangle +\frac{1}{2}(\sigma_{M}^{2}+\sigma_{I}^{2})\left\langle \phi'',\nu_{s}\right\rangle]ds+\sigma_{M}\int_{0}^{t}\left\langle \phi',\nu_{s}\right\rangle dW_{s},
\]
for all $t\in\left[0,T\right]$ and $\phi\in C^{\textrm{test}}$.
\end{prop}

\begin{proof} Take $\phi\in C^{\textrm{test}}$, then
It\=o's formula gives 
\begin{align}
\phi(X_{t\wedge\tau^{1}}^{1}) &=  \phi(X_{0}^{1})+\int_{0}^{t}[\mu\phi'(X_{s}^{1})+\frac{1}{2}(\sigma_{M}^{2}+\sigma_{I}^{2})\phi''(X_{s}^{1})]\mathbf{1}_{s<\tau^{1}}ds\label{eq:ItoForm1}\\
  & +\sigma_{M}\int_{0}^{t\wedge\tau^{1}}\phi'(X_{s}^{1})\mathbf{1}_{s<\tau^{1}}dW_{s}+\sigma_{I}\int_{0}^{t\wedge\tau^{1}}\phi'(X_{s}^{1})\mathbf{1}_{s<\tau^{1}}dW_{s}^{1}.\nonumber 
\end{align}
Since $\phi(0)=0$, we can write $\phi(X_{t\wedge\tau^{1}}^{1})=\phi(X_{t}^{1})\mathbf{1}_{t<\tau^{1}}$,
and therefore, by taking the conditional expectation over (\ref{eq:ItoForm1})
and noting that the $W^{1}$ integral then vanishes, we arrive at
the required result. \end{proof}

\begin{rem}[The choice of test functions] $C^{\mathrm{test}}$
is chosen as it is the largest space for which the previous proof
remains valid. The test functions' derivatives need not be controlled
at the origin, and we cannot take $C^{\mathrm{test}}=C_{0}^{\infty}(0,\infty)$
--- the space of compactly supported smooth functions on the half-line
--- because in section 4 it will be necessary to set $T_{\delta}\phi$
(see Definition \ref{HeatKernel}) into (\ref{eq:TheWeakProblem})
despite the fact that in general $\left(T_{\delta}\phi\right)'\left(0\right)$
and $\left(T_{\delta}\phi\right)''\left(0\right)$ are non-zero.
\end{rem}

\begin{rem}[Connection with filtering] If we treat $X^{1}$
as a signal process and $W$ an observation process, then the Zakai
equation (see for example \cite{Bain}) for the conditional distribution
of the signal given the observations is exactly (\ref{eq:TheWeakProblem}).
\end{rem}

\subsection{Properties}

Conditional on knowing $\mathcal{F}_{t}^{W}$, the process $X^{1}$
is a Brownian motion started at $X_{0}^{1}$ with drift $t\mapsto\mu t+\sigma_{M}W_{t}$,
and therefore we can write the conditional density of $X^{1}$ as
\[
P_{t}\left(x;x_{0}\right)=P_{t}\left(x;x_{0}\right)\left(\omega\right):=\frac{1}{\sqrt{2\pi\sigma_{I}^{2}t}}\exp\left\{ -\frac{\left(x-x_{0}-\mu t-\sigma_{M}W_{t}\right)^{2}}{2\sigma_{I}^{2}t}\right\} .
\]
From Definition \ref{ProbSolution}, it is a trivial fact that $\nu_{t}\left(S\right)\leq\mathbb{P}\left(\left.X_{t}^{1}\in S\right|\mathcal{F}_{t}^{W}\right)$
for any time $t\in\left[0,T\right]$ and measurable subset $S\subseteq\left(0,\infty\right)$,
and this simple observation allows us to show that $\nu$ has a density
process:

\begin{prop}[Existence of the density]
\label{ExistenceRad-Nik}There exists a full subset of $\Omega$ on which, for every $t\in\left[0,T\right]$,
$\nu_{t}$ has a density process $V_{t}$, and, for every $1\leq p\leq\infty$,
$V_{t}\in L^{p}\left(0,\infty\right)$ and $\left\Vert V_{t}\right\Vert _{p}\leq\left\Vert V_{0}\right\Vert _{p}$.
\end{prop} 

\begin{proof} 
From the above, there exists a full subset
of $\Omega$ on which 
\[
\nu_{t}\left(S\right)\leq\mathbb{P}\left(\left.X_{t}^{1}\in S\right|\mathcal{F}_{t}^{W}\right)=\int_{S}\int_{0}^{\infty}P_{t}\left(x;x_{0}\right)V_{0}\left(x_{0}\right)dx_{0}dx
\]
for all $t\in\left[0,T\right]$ and $S\subseteq\left(0,\infty\right)$
measurable. Therefore, if $S\subseteq\left(0,\infty\right)$ has zero
Lebesgue measure, then $\nu_{t}\left(S\right)=0$. It follows by the
Radon--Nikodym Theorem that $\nu_{t}$ has a density for every $t\in\left[0,T\right]$,
and that this density $V_{t}$ satisfies 
\begin{equation}
V_{t}\left(x\right)\leq\int_{0}^{\infty}P_{t}\left(x;x_{0}\right)V_{0}\left(x_{0}\right)dx_{0}\label{eq:DenCom}
\end{equation}
for almost all $x\in\left(0,\infty\right)$.

Bounding $V_{0}$ by $\left\Vert V_{0}\right\Vert _{\infty}$ in (\ref{eq:DenCom}),
we get $V_{t}\in L^{\infty}\left(0,\infty\right)$ with $\left\Vert V_{t}\right\Vert _{\infty}\leq\left\Vert V_{0}\right\Vert _{\infty}$.
Suppose $1\leq p<\infty$, then Holder's inequality applied to (\ref{eq:DenCom})
gives 
\[
V_{t}\left(x\right)^{p}\leq\int_{0}^{\infty}P_{t}\left(x;x_{0}\right)V_{0}\left(x_{0}\right)^{p}dx_{0}.
\]
Integrating over $x\in\left(0,\infty\right)$ and noting that $\int_{0}^{\infty}P_{t}\left(x;x_{0}\right)dx\leq1$
gives $\left\Vert V_{t}\right\Vert _{p}\leq\left\Vert V_{0}\right\Vert _{p}$.
\end{proof}

The most important feature of the probabilistic solution is the following
estimate:

\begin{prop}[The $\left(3+\beta\right)$-condition]
\label{TDCE2}Let $V_{0}\in L^{2}\left(0,\infty\right)$ be bounded and $\alpha$
be defined as in (\ref{eq:alpha}). Then for every $\beta\in\left(0,\pi/\alpha-1\right)\subseteq\left(0,1\right)$
there exists constants $B>0$ and $\gamma\in\left(0,1\right)$ such
that
\[
\mathbb{E}\left[\nu_{t}\left(0,\varepsilon\right)^{2}\right]\leq\frac{B}{t^{\gamma}}\varepsilon^{3+\beta},\quad\textrm{for every } \varepsilon>0 \textrm{ and } t\in\left[0,T\right].
\]
\end{prop} 

\noindent A proof is presented in Section 3.

In Lemma 3.5 of \cite{Bush}, this result is presented for an initial
density supported on a finite interval away from the boundary, and
that assumption allows the authors to drop the singular factor of
$t^{-\gamma}$. However, this is not the most natural initial condition;
for any $t>0$, $V_{t}$ takes positive values on all of $\left(0,\infty\right)$,
so if we stop the process at a time $t_{0}>0$ and then restart from
the density $V_{t_{0}}$, we still have the estimate on $\mathbb{E}\left[\nu_{t}\left(0,\varepsilon\right)^{2}\right]$
for $t>t_{0}$ despite $V_{t_{0}}$ not being supported on a finite
interval away from zero. In this more general setting, we can expect
a singular time factor to appear, since, for every time $t>0$, $V_{t}\left(x\right)\rightarrow0$
as $x\rightarrow0$ (Theorem 3.2 of \cite{KrylovHeat}) and so positive
values of $V_{0}$ close to zero must decay instantaneously as the
system evolves in time. Since $\gamma\in\left(0,1\right)$, $t^{-\gamma}$
is integrable over $\left[0,T\right]$ and therefore does not present
technical difficulties in later proofs of $L^{2}$-integrability. 

Using Proposition \ref{HittingProb2}, we prove Proposition \ref{TDCE2}
by considering the process $\mathbf{X}_{t}:=\left(X_{t}^{1},X_{t}^{2}\right)$,
which is a two-dimensional Brownian motion with components of correlation
$\rho$ --- recall (\ref{eq:alpha}). With $\alpha$ also as in (\ref{eq:alpha}),
the map $F_{\alpha}:\mathbb{R}^{2}\rightarrow\mathbb{R}^{2}$ defined
by 
\begin{equation}
F_{\alpha}:\mathbf{x}\mapsto\left(\begin{array}{cc}
\sqrt{1-\rho^{2}} & \rho\\
0 & 1
\end{array}\right)^{-1}\mathbf{x}\label{eq:AngleTransform}
\end{equation}
transforms $\mathbf{X}$ to a Brownian motion with uncorrelated components,
and maps the quadrant $\left(0,\infty\right)^{2}$ to the wedge 
\[
F_{\alpha}\left[\left(0,\infty\right)^{2}\right]=\left\{ \left(r,\theta\right):0<r<\infty,\,0<\theta<\alpha\right\} .
\]
We make use of explicit formulae in \cite{Iyengar} and \cite{Metzler}
to estimate the probability that $F_{\alpha}\left(\mathbf{X}_{t}\right)$
is in a small neighbourhood of the apex of the wedge and has not exited
the wedge, which corresponds to the event $X_{t}^{1},X_{t}^{2}\in\left(0,\varepsilon\right)$
and $t<\min\left(\tau^{1},\tau^{2}\right)$. The kernel smoothing
technique of Section 4, which was introduced in \cite{Bush} and adapted
from \cite{Kurtz}, relates this hitting probability to the behaviour
of the solution near the absorbing boundary.


\section{Proof of Proposition \ref{TDCE2}}

Let us write $\mathbb{P}_{x_{1},x_{2}}$ to indicate that $X_{0}^{1}=x_{1}$
and $X_{0}^{2}=x_{2}$, that is, we start the two-dimensional Brownian
motion, $\mathbf{X}=\left(X^{1},X^{2}\right)$, from $\left(x_{1},x_{2}\right)$.
From Proposition \ref{HittingProb2}, we know 
\begin{equation}
\mathbb{E}_{x_{1},x_{2}}\left[\nu_{t}\left(S_{1}\right)\nu_{t}\left(S_{2}\right)\right]=\mathbb{E}_{x_{1},x_{2}}\left[\mathbb{E}\left[\left.\mathbf{1}_{X_{t}^{1}\in S_{1},\, X_{t}^{2}\in S_{2}},t<\min\left(\tau^{1},\tau^{2}\right)\right|\mathcal{F}_{t}^{W}\right]\right],\label{eq:2dDensity}
\end{equation}
for every measurable $S_{1},S_{2}\subseteq\left(0,\infty\right)$.
A change of measure allows us to eliminate the drift term in the dynamics
of $\mathbf{X}$:

\begin{lem}[Eliminating drift]
\label{DriftKill} Let $\tilde{\mathbb{P}}$ be defined on $\left(\Omega,\mathcal{F}_{t}\right)$
by the Girsanov transformation 
\[
\frac{d\tilde{\mathbb{P}}}{d\mathbb{P}}=Z_{t}\left(\mu\right):=\exp\left\{ -\frac{\mu}{\sigma_{M}}W_{t}-\frac{\mu^{2}}{2\sigma_{M}^{2}}t\right\} ,
\]
then under $\tilde{\mathbb{P}}$ $\mathbf{X}$ is a two-dimensional
Brownian motion with zero drift. Let $\tilde{\mathbb{E}}$ denote
the expectation under $\tilde{\mathbb{P}}$. If there exists a constant
$C>0$ such that 
\[
\tilde{\mathbb{E}}\left[\nu_{t}\left(0,\varepsilon\right)^{2}\right]\leq\frac{C}{t^{\pi/2\alpha}}\varepsilon^{2+\pi/\alpha},\quad\mathrm{for\, all}\, t\in\left(0,T\right]\,\mathrm{and}\,\varepsilon>0,
\]
then Proposition \ref{TDCE2} holds.
\end{lem} 

\begin{proof} 
Following the method in Lemma 3.5 of \cite{Bush}, for any H\"older conjugates
$1<a,b<\infty$ we have 
\begin{align*}
\mathbb{E}\bigl[\nu_{t}\bigl(0,\varepsilon\bigr)^{2}\bigr]  
  & = \tilde{\mathbb{E}}\bigl[Z_{t}\bigl(\mu\bigr)\cdot\nu_{t}\bigl(0,\varepsilon\bigr)^{2}\bigr] \\
  & \leq  \tilde{\mathbb{E}}\bigl[Z_{t}\bigl(\mu\bigr)^{a}\bigr]^{\frac{1}{a}}\tilde{\mathbb{E}}\bigl[\nu_{t}\bigl(0,\varepsilon\bigr)^{2b}\bigr]^{\frac{1}{b}}\\
  & \leq  \tilde{\mathbb{E}}\bigl[Z_{t}\bigl(\mu\bigr)^{a}\bigr]^{\frac{1}{a}}\tilde{\mathbb{E}}\bigl[\nu_{t}\bigl(0,\varepsilon\bigr)^{2}\bigr]^{\frac{1}{b}}\\
 & \leq \tilde{\mathbb{E}}\left[Z_{t}\left(\mu\right)^{a}\right]^{\frac{1}{a}}\bigl\{ Ct^{-\pi/2\alpha}\varepsilon^{2+\pi/\alpha}\bigr\} ^{\frac{1}{b}}
  =  \mathbb{E}\left[Z_{t}\left(\mu\right)^{a-1}\right]^{\frac{1}{a}}\bigl\{ Ct^{-\pi/2\alpha}\varepsilon^{2+\pi/\alpha}\bigr\} ^{\frac{1}{b}},
\end{align*}
where one should note that $\nu_{t}\left(0,\varepsilon\right)\leq1$.
Hence 
\[
\mathbb{E}\left[\nu_{t}\left(0,\varepsilon\right)^{2}\right]\leq\exp\left(\frac{\left(a-1\right)\left(a-2\right)\mu^{2}}{2\sigma_{M}^{2}a}T\right)\left\{ \frac{C}{t^{\pi/2\alpha}}\varepsilon^{2+\pi/\alpha}\right\} ^{\frac{1}{b}},
\]
and, since $\alpha\in\left(\pi/2,\pi\right)$, we can therefore choose
$b$ sufficiently close to 1 so that the result holds. 
\end{proof}

\begin{rem} We cannot take $b=1$ above, and so, except
in the case when $\mu=0$, this change of measure removes the possibility
of considering the borderline case $\beta=\pi/\alpha-1$. Although
this may seem to weaken the potential regularity available near the
absorbing boundary, in proving Theorem \ref{Reg} we make use
of Lemma \ref{RemainderTerm} which requires $\beta$ to be strictly
smaller than the maximal value, $\pi/\alpha-1$, and hence it is unclear
whether we could analyse the case $\beta=\pi/\alpha-1$ regardless
of this difficulty caused by eliminating the drift.
\end{rem}

The remainder of Section 3 is devoted to showing that the hypothesis
of Lemma \ref{DriftKill} holds. Without loss of generality, from
this point on we take $\mu=0$ and $\tilde{\mathbb{P}}=\mathbb{P}$.

Conditioning on the start point of $\mathbf{X}$, we have 
\[
\mathbb{E}\left[\nu_{t}\left(0,\varepsilon\right)^{2}\right]=\int_{0}^{\infty}\int_{0}^{\infty}\mathbb{E}_{x_{1},x_{2}}\left[\mathbb{E}\left[\left.\mathbf{1}_{X_{t}^{1}\in\left(0,\varepsilon\right),\, X_{t}^{2}\in\left(0,\varepsilon\right)}\right|\mathcal{F}_{t}^{W}\right]\right]V_{0}\left(x_{1}\right)V_{0}\left(x_{2}\right)dx_{1}dx_{2}.
\]
Under the transformation $F_{\alpha}$, the region $\left(0,\varepsilon\right)^{2}$
is mapped \emph{into} a region of the form 
\[
\left\{ \left(r,\theta\right):0<r<c\varepsilon,0<\theta<\alpha\right\} ,
\]
for some numerical constant $c$. The estimate in Proposition \ref{TDCE2}
is unchanged under the map $\varepsilon\mapsto c\varepsilon$ (modulo
multiplicative constants), and so it is no loss of generality to assume
$c=1$. Using that $V_{0}$ is bounded, that the Jacobian of the transform
$F_{\alpha}:\mathbb{R}^{2}\rightarrow\mathbb{R}^{2}$ is a numerical
constant depending only on $\rho$, and formula (8) of \cite{Iyengar}
(corroborated in \cite{Metzler}), we arrive at 
\begin{align*}
\mathbb{E} \! \left[\nu_{t}\left(0,\varepsilon\right)^{2}\right] 
  & \leq C \!\! \int_{0}^{\infty} \!\!\!\! \int_{0}^{\alpha} \!\!\!\! \int_{0}^{\alpha} \!\!\!\! \int_{0}^{\varepsilon}\frac{rr_{0}}{t}e^{-\frac{r^{2}+r_{0}^{2}}{2t}}\sum_{n=1}^{\infty}\sin\!\left(\frac{n\pi\theta}{\alpha}\right)\sin\!\left(\frac{n\pi\theta_{0}}{\alpha}\right)I_{\frac{n \pi}{\alpha}}\!\left(\frac{rr_{0}}{2t}\right)\!drd\theta dr_{0}d\theta_{0} \\
  & =  C\!\!\int_{0}^{\varepsilon}\frac{r^{\pi/\alpha+1}}{t^{\pi/2\alpha}}e^{-\frac{r^{2}}{2t}}\sum_{n,m=0}^{\infty}\frac{\Gamma\left(m+\frac{\left(2n+1\right)\pi}{2\alpha}+1\right)}{\left(2n+1\right)^{2}m!\Gamma\left(m+\frac{\left(2n+1\right)\pi}{\alpha}+1\right)}\left(\frac{r}{2t}\right)^{m+n\pi/\alpha}dr,
\end{align*}
for $C>0$ some numerical constant. As $\sup_{x\geq0}\left\{ x^{c}e^{-x}\right\} =c^{c}e^{-c}$,
we have 
\begin{equation}
\mathbb{E}\left[\nu_{t}\left(0,\varepsilon\right)^{2}\right]\leq C\int_{0}^{\varepsilon}\frac{r^{\pi/\alpha+1}}{t^{\pi/2\alpha}}dr\sum_{n,m=0}^{\infty}\frac{\Gamma\left(m+\frac{\left(2n+1\right)\pi}{2\alpha}+1\right)\left(m+\frac{n\pi}{\alpha}\right)^{m+\frac{n\pi}{\alpha}}e^{-\left(m+\frac{n\pi}{\alpha}\right)}}{\left(2n+1\right)^{2}m!\Gamma\left(m+\frac{\left(2n+1\right)\pi}{\alpha}+1\right)},\label{eq:DSum}
\end{equation}
and therefore we are done provided the double sum on the right-hand
is finite. Using Stirling's approximation, the inequality 
\[
\frac{\Gamma\left(a+c+1\right)}{\Gamma\left(b+c+1\right)}\leq\frac{\Gamma\left(a+1\right)}{\Gamma\left(b+1\right)}
\]
for $a,b,c>0$, and the asymptotic behaviour $\Gamma\left(m+\pi/2\alpha+1\right)\sim m!m^{\pi/2\alpha}$
for large $m$, we have that, for large $m$ and $m+n$, the summand
in (\ref{eq:DSum}) is dominated by a constant multiple of 
\[
\frac{1}{\left(2n+1\right)^{2}m^{\pi/2\alpha+1/2}}.
\]
Since $\pi/2\alpha+1/2>1$, this expression is summable over $m$
and $n$, and so we have the result. \qed


\section{Proof of Theorem \ref{Reg}}

We begin by proving the regularity result for the probabilistic solution.
Let us write 
\begin{equation}
p_{\delta}\left(x\right)=\frac{1}{\sqrt{2\pi\delta}}\exp\left\{ -\frac{x^{2}}{2\delta}\right\}, \quad \textrm{for } z\in\mathbb{R} \textrm{ and } \delta>0.\label{eq:TransitionDenisty}
\end{equation}

\begin{defn}[Heat kernels]
\label{HeatKernel}The \emph{absorbing} and \emph{reflected heat kernels} are defined
to be
\[
G_{\delta}\left(x,y\right):=p_{\delta}\left(x-y\right)-p_{\delta}\left(x+y\right)\quad\mathrm{and}\quad G_{\delta}^{r}\left(x,y\right):=p_{\delta}\left(x-y\right)+p_{\delta}\left(x+y\right)
\]
 for $x,y>0$ and $\delta>0$. The \emph{smoothed density} and \emph{reflected
smoothed density} are then defined to be 
\[
T_{t,\delta}\left(x\right):=\int_{0}^{\infty}G_{\delta}\left(x,y\right)V_{t}\left(y\right)dy\quad\mathrm{and}\quad T_{t,\delta}^{r}\left(x\right):=\int_{0}^{\infty}G_{\delta}^{r}\left(x,y\right)V_{t}\left(y\right)dy
\]
for $t\in\left[0,T\right]$, $x>0$, and $\delta>0$.
\end{defn}

Take any smooth compactly supported function, $\phi\in C_{0}^{\infty}(0,\infty)$,
then, for every $\delta>0$, $x\mapsto T_{\delta}\phi(x)$
(that is the function obtained by convolving the heat kernel, $G_{\delta}$,
with $\phi$) is in $C^{\textrm{test}}$, so putting $T_{\delta}\phi$
and $\nu$ into (\ref{eq:TheWeakProblem}) gives the following strong
equation for the evolution of the smoothed density:

\begin{prop}[Strong smoothed evolution equation]
\label{SEE}For $\delta>0$, $t\in\left[0,T\right]$, and $x>0$ 
\begin{align}
T_{t,\delta}\left(x\right)
   & = T_{0,\delta}\left(x\right)+\int_{0}^{t}\bigl[-\mu T_{s,\delta}'\left(x\right)+\frac{\sigma_{M}^{2}+\sigma_{I}^{2}}{2}T_{s,\delta}''\left(x\right)\bigr]ds-\sigma_{M}\int_{0}^{t}T_{s,\delta}'\left(x\right)dW_{s}\label{eq:SEEwR}\\
 & +\mu\int_{0}^{t}R_{s,\delta}'\left(x\right)ds+\sigma_{M}\int_{0}^{t}R_{s,\delta}'\left(x\right)dW_{s},\nonumber 
\end{align}
where $R_{t,\delta}$ is the \emph{remainder function}: 
\[
R_{t,\delta}\left(x\right):=T_{t,\delta}\left(x\right)-T_{t,\delta}^{r}\left(x\right)=2\int_{0}^{\infty}p_{\delta}\left(x+y\right)V_{t}\left(y\right)dy.
\]
\end{prop}

\begin{proof} This is simply a matter of observing that
\begin{eqnarray*}
\bigl\langle \phi',T_{t,\delta}\bigr\rangle  & = & -\bigl\langle \phi,\left(T_{t,\delta}^{r}\right)'\bigr\rangle \\
\bigl\langle \phi'',T_{t,\delta}\bigr\rangle  & = & \bigl\langle \phi,T_{t,\delta}''\bigr\rangle .
\end{eqnarray*}
One also needs to apply the stochastic Fubini theorem \cite{Bichteler}
to switch integration in the space dimension with integration with
respect to $W$ (and likewise with respect to time). To this end,
it is enough to note that the tails of $T_{t,\delta}$ and $T_{t,\delta}^{r}$
and their derivatives decay exponentially, and so the conditions of
the (stochastic) Fubini theorem are met.
\end{proof}

Differentiating equation (\ref{eq:SEEwR}) $n$ times (or integrating
once for the case $n=-1$ --- the result is no different) and applying
It\=o's formula to $\bigl(T_{t,\delta}^{\left(n\right)}\bigr)^{2}$
gives 
\begin{align}
d\bigl(T_{t,\delta}^{\left(n\right)}\bigr)^{2} 
  & = -2\mu T_{t,\delta}^{\left(n\right)}T_{t,\delta}^{\left(n+1\right)}dt+\left(\sigma_{M}^{2}+\sigma_{I}^{2}\right)T_{t,\delta}^{\left(n\right)}T_{t,\delta}^{\left(n+2\right)}dt+\sigma_{M}^{2}\bigl(T_{t,\delta}^{\left(n+1\right)}\bigr)^{2}dt \nonumber \\
  & +\sigma_{M}^{2}\bigl(R_{t,\delta}^{\left(n+1\right)}\bigr)^{2}dt+2\mu T_{t,\delta}^{\left(n\right)}R_{t,\delta}^{\left(n+1\right)}dt\label{eq:Inter}\\
 & -2\sigma_{M}T_{t,\delta}^{\left(n\right)}T_{t,\delta}^{\left(n+1\right)}dW_{t}+2\sigma_{M}T_{t,\delta}^{\left(n\right)}R_{t,\delta}^{\left(n+1\right)}dW_{t}.\nonumber 
\end{align}
From this equation, we shall proceed to a proof by induction in Section
4.2. Before continuing, we introduce several technical lemmas in the
next section. 

\subsection*{Some technical lemmas}

\begin{lem}
\label{E1}There exists a full subset of $\Omega$ on which $x^{-c}T_{t,\delta}\left(x\right)\rightarrow0$
as $x\rightarrow0$ for every $c\in\left(0,1\right)$, $\delta>0$,
and $t\in\left(0,T\right]$.
\end{lem}

\begin{proof}
This is straightforward
since, for $x<1$, we have 
\begin{align*}
0 \leq T_{t,\delta}\left(x\right) & 
  = \frac{1}{\sqrt{2\pi\delta}}\int_{0}^{\infty}\Bigl[\exp\Bigl\{ -\frac{\left(x-y\right)^{2}}{2\delta}\Bigr\} -\exp\Bigl\{ -\frac{\left(x+y\right)^{2}}{2\delta}\Bigr\} \Bigr]V_{t}\left(y\right)dy\\
 & \leq \frac{2\left\Vert V_{0}\right\Vert _{\infty}x}{\delta\sqrt{2\pi\delta}}\int_{0}^{\infty}y\exp\Bigl\{ -\frac{\left(y-1\right)^{2}}{2\delta}\Bigr\} dy = D_{\delta}x,
\end{align*}
for a numerical constant $D_{\delta}>0$. Note that we used the estimate
$1-e^{-z}\leq z$.
\end{proof}

\begin{lem}
\label{liminfLemma}Let $n\geq0$ be an integer and suppose that for $k=0,1,\ldots,n$
\[
\liminf_{\delta\rightarrow0}\mathbb{E}\int_{0}^{T}\bigl\Vert w_{k-1-\beta/2}\left(T_{t,\delta}\right)^{\left(k\right)}\bigr\Vert _{2}^{2}dt < \infty.
\]
Then, for almost every $\left(\omega,t\right)\in\Omega\times\left[0,T\right]$,
the density $V_{t}$ is $n$-times weakly differentiable with 
\begin{equation}
\mathbb{E}\int_{0}^{T}\bigl\Vert w_{k-1-\beta/2}V_{t}^{\left(k\right)}\bigr\Vert_{2}^{2}dt<\infty,\quad\textrm{for } k=0,1,\ldots,n.\label{eq:WeightedNorm}
\end{equation}
\end{lem}

\begin{proof} Let $\left\langle \cdot,\cdot\right\rangle $
denote the usual $L^{2}(0,\infty)$ inner product: $\left\langle f,g\right\rangle =\int_{0}^{\infty}f\left(x\right)g\left(x\right)dx$.
Let $\left\{ \phi_{i}\right\} _{i=1}^{\infty}$ be a basis of $L^{2}(0,\infty)$
for which each $\phi_{i}$ is smooth and bounded. 

Firstly, for any $\phi\in C^{\infty}(0,\infty)$ and $m>k\geq1$,
$\lim_{x\downarrow0}\left(w_{m-\beta/2}\phi\left(x\right)\right)^{\left(k\right)}=0$,
since $\partial_{x}^{k}w_{m-\beta/2}$ is a linear combination of
$w_{i-\beta/2}$'s and $w_{i-\beta/2}\rightarrow0$ for $m-k\leq i\leq m$.
So for $k\geq1$, integrating by parts gives 
\[
\bigl\langle \left(w_{k-1-\beta/2}\phi\right)^{\left(k\right)},T_{t,\delta}\bigr\rangle  = \left(-1\right)^{k}\bigl\langle \phi,w_{k-1-\beta/2}T_{t,\delta}^{\left(k\right)}\bigr\rangle -\lim_{x\downarrow0}\bigl\{ T_{t,\delta}\left(x\right)\left(w_{k-1-\beta/2}\phi\right)^{\left(k-1\right)}\bigr\} .
\]
Applying Leibniz\textquoteright{}s rule to $\left(w_{k-1-\beta/2}\phi\right)^{\left(k-1\right)}$,
we see that 
\[
\lim_{x\downarrow0}\bigl\{ T_{t,\delta}\left(x\right)\left(w_{k-1-\beta/2}\phi\right)^{\left(k-1\right)}\left(x\right)\bigr\} =C_{k}\lim_{x\downarrow0}\left\{ T_{t,\delta}\left(x\right)w_{-\beta/2}\left(x\right)\phi\left(x\right)\right\} ,
\]
for some numerical constant $C_{k}$, and that the right-hand side
vanishes by Lemma \ref{E1}, because $w_{-\beta/2}\left(x\right)=\mathcal{O}(x^{-\beta/2})$
for small $x$. Hence we conclude 
\[
\bigl\langle \left(w_{k-1-\beta/2}\phi\right)^{\left(k\right)},T_{t,\delta}\bigr\rangle =\left(-1\right)^{k}\bigl\langle \phi,w_{k-1-\beta/2}T_{t,\delta}^{\left(k\right)}\bigr\rangle .
\]
The remaining case, $k=0$, is trivial. 

By the above results and Fatou's Lemma, we have 
\begin{align*}
\mathbb{E}\int_{0}^{T}\sum_{i=1}^{\infty}\bigl|\bigl\langle \left(w_{k-1-\beta/2}\phi_{i}\right)^{\left(k\right)},V_{t}\bigr\rangle \bigr|^{2}dt 
  & = \mathbb{E} \int_{0}^{T}\sum_{i=1}^{\infty}\lim_{\delta\rightarrow0}\bigl|\bigl\langle \left(w_{k-1-\beta/2}\phi_{i}\right)^{\left(k\right)},T_{t,\delta}\bigr\rangle \bigr|^{2}dt \\
  & \leq \liminf_{\delta\rightarrow0}\mathbb{E} \int_{0}^{T}\sum_{i=1}^{\infty}\bigl|\bigl\langle \phi_{i},w_{k-1-\beta/2}T_{t,\delta}^{\left(k\right)}\bigr\rangle \bigr|^{2}dt \\
 & = \liminf_{\delta\rightarrow0}\mathbb{E}\int_{0}^{T}\bigl\Vert w_{k-1-\beta/2}T_{t,\delta}^{\left(k\right)}\bigr\Vert _{2}^{2}dt < \infty,
\end{align*}
for $k=0,1,\ldots n$, and therefore we can define the processes 
\[
W_{t,k}\left(x\right):=\sum_{i=1}^{\infty}\bigl\langle \left(w_{k-1-\beta/2}\phi_{i}\right)^{\left(k\right)},V_{t}\bigr\rangle \phi_{i}\left(x\right)
\]
which satisfy 
\[
\mathbb{E}\int_{0}^{T}\left\Vert W_{t,k}\right\Vert _{2}^{2}dt <\infty.
\]

For $k=0$, let $f_{0}=w_{-1-\beta/2}$ and take any $\psi\in C_{0}^{\infty}\left(0,\infty\right)$,
then we calculate 
\[
\bigl\langle \psi,\frac{1}{f_{0}}W_{t,0}\bigr\rangle 
  = \sum_{i=1}^{\infty}\bigl\langle f_{0}\phi_{i},V_{t}\bigr\rangle \bigl\langle \frac{\psi}{f_{0}},\phi_{i}\bigr\rangle 
  = \Bigl\langle f_{0}\sum_{i=1}^{\infty}\bigl\langle \frac{\psi}{f_{0}},\phi_{i}\bigr\rangle _{2}\phi_{i},V_{t}\Bigr\rangle 
  = \bigl\langle \psi,V_{t}\bigr\rangle,
\]
and so $V_{t}=W_{t,0}/f_{0}$. Hence we have demonstrated (\ref{eq:WeightedNorm})
for $k=0$. 

Take $k>0$, and set $f_{k}=w_{k-1-\beta/2}$. Using $V_{t}=W_{t,0}/f_{0}$
we can compute 
\begin{align*}
\bigl\langle \psi,\frac{1}{f_{k}}W_{t,k}\bigr\rangle  
  & = \sum_{i=1}^{\infty}\bigl\langle \left(f_{k}\phi_{i}\right)^{\left(k\right)},V_{t}\bigr\rangle \bigl\langle \frac{\psi}{f_{k}},\phi_{i}\bigr\rangle \\
 & = \sum_{i,j=1}^{\infty}\bigl\langle \left(f_{k}\phi_{i}\right)^{\left(k\right)},\frac{\phi_{j}}{f_{0}}\bigr\rangle \bigl\langle f_{0}\phi_{j},V_{t}\bigr\rangle \bigl\langle \frac{\psi}{f_{k}},\phi_{i}\bigr\rangle \\
 & = \sum_{i,j,m=1}^{\infty}\bigl\langle \left(f_{k}\phi_{i}\right)^{\left(k\right)},\phi_{m}\bigr\rangle \bigl\langle \frac{\phi_{j}}{f_{0}},\phi_{m}\bigr\rangle \bigl\langle f_{0}\phi_{j},V_{t}\bigr\rangle \bigl\langle \frac{\psi}{f_{k}},\phi_{i}\bigr\rangle \\
 & = \sum_{i,j,m=1}^{\infty}\left(-1\right)^{k}\bigl\langle \phi_{i},f_{k}\phi_{m}^{\left(k\right)}\bigr\rangle \bigl\langle \frac{\phi_{j}}{f_{0}},\phi_{m}\bigr\rangle \bigl\langle f_{0}\phi_{j},V_{t}\bigr\rangle \bigl\langle \frac{\psi}{f_{k}},\phi_{i}\bigr\rangle \\
 & = \sum_{j,m=1}^{\infty}\left(-1\right)^{k}\bigl\langle \psi,\phi_{m}^{\left(k\right)}\bigr\rangle \bigl\langle \frac{\phi_{j}}{f_{0}},\phi_{m}\bigr\rangle \bigl\langle f_{0}\phi_{j},V_{t}\bigr\rangle \\
 & = \sum_{j=1}^{\infty}\bigl\langle \psi^{\left(k\right)},\frac{\phi_{j}}{f_{0}}\bigr\rangle \bigl\langle f_{0}\phi_{j},V_{t}\bigr\rangle 
 = \bigl\langle \psi^{\left(k\right)},V_{t}\bigr\rangle,
\end{align*}
hence $\left(-1\right)^{k}W_{t,k}/f_{k}$ is the $k^{\textrm{th}}$
weak derivative of $V_{t}$, and this completes the proof.
\end{proof}

\begin{lem}
\label{IBP}There exists a full subset of $\Omega$ on which
\[
\int_{0}^{\infty}w_{n-\beta/2}^{2}T_{t,\delta}^{\left(n\right)}(x)T_{t,\delta}^{\left(n+1\right)}(x)dx=\bigl\Vert w_{n-\beta/2}T_{t,\delta}^{\left(n\right)}\bigr\Vert _{2}^{2}-\bigl(n-\frac{\beta}{2}\bigr)\bigl\Vert w_{n-1/2-\beta/2}T_{t,\delta}^{\left(n\right)}\bigr\Vert _{2}^{2}
\]
and
\begin{align*}
\int_{0}^{\infty}w_{n-\beta/2}^{2}T_{t,\delta}^{\left(n\right)}(x) T_{t,\delta}^{\left(n+2\right)}(x) dx 
  & = 2\bigl\Vert w_{n-\beta/2}T_{t,\delta}^{\left(n\right)}\bigr\Vert _{2}^{2}-\bigl\Vert w_{n-\beta/2}T_{t,\delta}^{\left(n+1\right)}\bigr\Vert _{2}^{2}\\
  & \quad -4\bigl(n-\frac{\beta}{2}\bigr)\bigl\Vert w_{n-1/2-\beta/2}T_{t,\delta}^{\left(n\right)}\bigr\Vert _{2}^{2}\\
  & \quad +2\bigl(n-\frac{\beta}{2}\bigr)\bigl(n-\frac{\beta+1}{2}\bigr)\bigl\Vert w_{n-1-\beta/2}T_{t,\delta}^{\left(n\right)}\bigr\Vert_{2}^{2},
\end{align*}
for every $n\geq-1$, $\delta>0$, and $t\in\left[0,T\right]$.
\end{lem}

\begin{proof} The result is an integration-by-parts exercise provided
the following limits vanish: 
\[
\lim_{x\downarrow0}\left\{ w_{n-1/2-\beta/2}\left(x\right)T_{t,\delta}^{\left(n\right)}\left(x\right)\right\} 
\]
\[
\lim_{x\downarrow0}\left\{ w_{n-\beta/2}^{2}\left(x\right)T_{t,\delta}^{\left(n\right)}\left(x\right)T_{t,\delta}^{\left(n+1\right)}\left(x\right)\right\} 
\]
for $n\geq-1$. If $n\geq1$, then $w_{n-1/2-\beta/2},w_{n-\beta/2}\rightarrow0$
as $x\rightarrow0$, and the result is immediate. In the case $n=0$,
Lemma \ref{E1} gives $T_{t,\delta}\left(x\right)=\mathcal{O}\left(x\right)$,
which is sufficient for the limit to be zero. This also implies $T_{t,\delta}^{\left(-1\right)}\left(x\right)=\mathcal{O}\left(x^{2}\right)$,
so the result holds for $n=-1$ too.
\end{proof}

\begin{lem}
\label{DiffTransition}With $p_{\delta}$ defined as in (\ref{eq:TransitionDenisty}), for
every integer $n\geq0$ there exist constants $c_{i}^{n}$ for $n/2\leq i\leq n$
such that 
\[
\frac{\partial^{n}p_{\delta}}{\partial x^{n}}\left(x\right)=\sum_{n/2\leq i\leq n}c_{i}^{n}x^{2i-n}\delta^{-i}p_{\delta}\left(x\right),\quad\mathrm{for\, every}\,\delta>0\,\mathrm{and\,}x>0.
\]
\end{lem}

\begin{proof} This is a simple inductive argument. \end{proof}

\begin{lem}
\label{InitialDensity} For $n=-1,0,1,\ldots,N$ and $\beta<\pi/\alpha-1$ 
\[
\liminf_{\delta\rightarrow0}\bigl\Vert w_{n-\beta/2}T_{0,\delta}^{\left(n\right)}\bigr\Vert _{2}^{2}<\infty.
\]
\end{lem}

\begin{proof} The cases $n=0$ and $n=-1$ are simple.
Firstly, observe from Definition \ref{HeatKernel} that $\left\Vert T_{0,\delta}\right\Vert _{\infty}\leq\left\Vert V_{0}\right\Vert _{\infty}$,
so for all $\delta>0$ 
\[
\left\Vert w_{-\beta/2}T_{0,\delta}\right\Vert _{2}^{2}\leq\left\Vert V_{0}\right\Vert _{\infty}\int_{0}^{\infty}x^{-\beta}e^{-2x}dx<\infty,
\]
and therefore we have the result for $n=0$. The case $n=-1$ follows
by noting that $0\leq T_{0,\delta}^{-1}\left(x\right)\leq\left\Vert V_{0}\right\Vert _{\infty}x$
and applying the same argument. 

Now fix $n\geq1$. Begin by noting that $\partial_{x}^{n}G_{\delta}\left(x,y\right)=\partial_{y}^{n}G_{\delta}^{r^{n}}\left(x,y\right)$,
where $G^{r^{n}}=G$ if $n$ is even and $G^{r^{n}}=G^{r}$ if $n$
is odd. Splitting the range of integration and integrating by parts
in the definition of $T_{0,\delta}$ gives 
\begin{align}
T_{0,\delta}^{\left(n\right)}\left(x\right) & = \sum_{i=0}^{n-1}\left(-1\right)^{n-1-i}\partial_{y}^{n-1-i}G_{\delta}^{r^{n}}\left(x,x/2\right)V_{0}^{\left(i\right)}\left(x/2\right)+\int_{x/2}^{\infty}G_{\delta}^{r^{n}}\left(x,y\right)V_{0}^{\left(n\right)}\left(y\right)dy\label{eq:Rewrite1} \nonumber  \\
 & \qquad +\int_{0}^{x/2}\partial_{x}^{n}G_{\delta}\left(x,y\right)V_{0}\left(y\right)dy.
\end{align}
We proceed by considering the weighted $L^{2}$-norm of these three
components. 

Firstly by Lemma \ref{DiffTransition}, there exists a numerical constant
$D>0$ such that 
\begin{multline*}
 \bigl|\partial_{y}^{n-1-i}G_{\delta}^{r^{n}}\left(x,x/2\right)V_{0}^{\left(i\right)}\left(x/2\right)w_{n-\beta/2}\left(x\right)\bigr|\\
  \leq D \!\!\!\! \sum_{\frac{n-1-i}{2}\leq j\leq n-1-i}x^{2j-n+1+i}\delta^{-j-1/2}e^{-\frac{x^{2}}{8\delta}}\bigl|V_{0}^{\left(i\right)}\left(x/2\right)\bigr|w_{i-\beta/2}\left(x\right)x^{n-i},
\end{multline*}
therefore integrating over $x\in\left(0,\infty\right)$ and applying
the Cauchy--Schwarz inequality gives 
\begin{eqnarray}
 &  & \bigl\Vert \partial_{y}^{n-1-i}G_{\delta}^{r^{n}}\left(\cdot,\cdot/2\right)V_{0}^{\left(i\right)}\left(\cdot/2\right)w_{n-\beta/2}\bigr\Vert _{2}^{2}\nonumber \\
 &  & \qquad\qquad\qquad \leq2D\sum_{\frac{n-1-i}{2}\leq j\leq n-1-i}\int_{0}^{\infty}x^{4j+2}\delta^{-2j-1}e^{-\frac{x^{2}}{4\delta}}dx\bigl\Vert V_{0}^{\left(i\right)}\left(\cdot/2\right)w_{i-\beta/2}\bigr\Vert _{2}^{2}\nonumber \\
 &  & \qquad\qquad\qquad\leq D'\sum_{\frac{n-1-i}{2}\leq j\leq n-1-i}\delta^{1/2}2^{2i-\beta}\bigl\Vert V_{0}^{\left(i\right)}w_{i-\beta/2}\bigr\Vert _{2}^{2}\rightarrow0\label{eq:Rewrite2}
\end{eqnarray}
as $\delta\rightarrow0$, where $D'>0$ is a further numerical constant. 

For the second term, we apply Cauchy--Schwarz to obtain 
\begin{eqnarray*}
 &  & \left\{ \int_{x/2}^{\infty}G_{\delta}^{r^{n}}(\cdot,y)V_{0}^{\left(n\right)}\left(y\right)dy\right\} ^{2}\\
 &  & \qquad \leq\left\{ \int_{0}^{\infty} \!\! G_{\delta}^{r^{n}}(\cdot,y)w_{n-\beta/2}\left(y\right)^{2}V_{0}^{\left(n\right)}\left(y\right)^{2}dy\right\} \left\{ \int_{x/2}^{\infty} \!\! G_{\delta}^{r^{n}}(\cdot,y)y^{-2n+\beta}e^{2y}dy\right\} \\
 &  & \qquad \leq4\cdot2^{2n-\beta}x^{-2n+\beta}\!\left\{ \int_{0}^{\infty} \!\!\! p_{\delta}\!\left(x-y\right)w_{n-\beta/2}\!\left(y\right)^{2}\!V_{0}^{\left(n\right)}\!\!\left(y\right)^{2}\!dy\right\} \! \left\{ \int_{x/2}^{\infty}\!\!p_{\delta}\left(x-y\right)e^{2y}dy\!\right\}\!.
\end{eqnarray*}
For $\delta<1$ we have $\int_{x/2}^{\infty}p_{\delta}\left(x-y\right)e^{2y}dy\leq D'e^{2x}$,
for some numerical constant $D'>0$, therefore 
\begin{align}
\Bigl\Vert w_{n-\beta/2}\int_{\cdot/2}^{\infty}\!\!G_{\delta}^{r^{n}}\!\left(\cdot,y\right)V_{0}^{\left(n\right)}\!\left(y\right)\!dy\Bigr\Vert _{2}^{2} 
  & \leq 2^{2n-\beta+2}D'\!\!\int_{0}^{\infty} \!\!\!\! \int_{0}^{\infty} \!\!\! p_{\delta}\!\left(x-y\right)\!dxw_{n-\beta/2}\!\left(y\right)^{2}V_{0}^{\left(n\right)}\!\left(y\right)^{2}\!dy\nonumber \\
  & \leq 2^{2n-\beta+2}D'\Bigl\Vert w_{n-\beta/2}V_{0}^{\left(n\right)}\Bigr\Vert _{2}^{2}.\label{eq:Rewrite3}
\end{align}
 
Finally, using the fact that $V_{0}$ is bounded and Lemma \ref{DiffTransition},
there exists a numerical constant $D''>0$ such that 
\[
\Bigl|\int_{0}^{x/2}\partial_{x}^{n}G_{\delta}\left(x,y\right)V_{0}\left(y\right)dy\Bigr| 
  \leq D''\left\Vert V_{0}\right\Vert _{\infty}x\sum_{n/2\leq i\leq n}x^{2i-n}\delta^{-i-1/2}e^{-\frac{x^{2}}{8\delta}},
\]
and therefore 
\begin{align}
\Bigl\Vert w_{n-\beta/2}\int_{0}^{\cdot/2}\partial_{x}^{n}G_{\delta}\left(\cdot,y\right)V_{0}\left(y\right)dy\Bigr\Vert _{2}^{2} 
  & \leq 2\left(D''\right)^{2}\left\Vert V_{0}\right\Vert _{\infty}^{2} \!\!\! \sum_{n/2\leq i\leq n}\int_{0}^{\infty}x^{4i-2n+2}\delta^{-2i-1}e^{-\frac{x^{2}}{8\delta}}dx\nonumber \\
  & \leq D'''\delta^{\left(1-\beta\right)/2}\rightarrow0\label{eq:Rewrite5}
\end{align}
as $\delta\rightarrow0$, since $\beta<1$. Here $D'''>0$ is a further
numerical constant. 

By using (\ref{eq:Rewrite1}), results (\ref{eq:Rewrite2}), (\ref{eq:Rewrite3}),
and (\ref{eq:Rewrite5}) complete the proof.
\end{proof}

\begin{lem}
\label{RemainderTerm}For $n\geq0$ and $\beta<\pi/\alpha-1$ 
\[
\liminf_{\delta\rightarrow0}\mathbb{E}\int_{0}^{t}\bigl\Vert w_{n-\beta/2}R_{s,\delta}^{\left(n+1\right)}\bigr\Vert _{2}^{2}ds=0.
\]
\end{lem}

\begin{proof}
Interchanging differentiation and integration and applying Lemma \ref{DiffTransition} gives 
\begin{align*}
\Bigl|R_{s,\delta}^{\left(n+1\right)}\left(x\right)\Bigr| & = \Bigl|2\int_{0}^{\infty}\frac{\partial^{n+1}p_{\delta}}{\partial x^{n+1}}\left(x+y\right)V_{s}\left(y\right)dy\Bigr|\\
 & \leq D\sum_{\frac{n+1}{2}\leq i\leq n+1}\int_{0}^{\infty}\left(x+y\right)^{2i-n-1}\delta^{-i-1/2}e^{-\frac{x^{2}}{2\delta}}e^{-\frac{y^{2}}{2\delta}}V_{s}\left(y\right)dy,
\end{align*}
with $D>0$ a numerical constant. Let $\eta\in\left(0,1\right)$ and
split the region of integration at $y=\delta^{\frac{1}{2}\left(1-\eta\right)}$,
then for $\delta<1$ we have $\delta^{-i}\leq\delta^{-n-1}$ and $e^{-\frac{y^{2}}{2\delta}}\leq e^{-\frac{y^{2}}{4}}e^{-1/4\delta^{\eta}}$
on $y>\delta^{\frac{1}{2}\left(1-\eta\right)}$, and so we have 
\begin{align}
 \Bigl|w_{n-\beta/2}(x)R_{s,\delta}^{\left(n+1\right)}(x)\Bigr|
 & \leq\Bigl\{ D\sum_{\frac{n+1}{2}\leq i\leq n+1}\int_{\delta^{\left(1-\eta\right)/2}}^{\infty}\left(x+y\right)^{2i-n-1}e^{-\frac{x^{2}}{2}}e^{-\frac{y^{2}}{4}}dy\Bigr\} \nonumber
 \\
 & \qquad \times  \left\Vert V_{0}\right\Vert _{\infty}\delta^{-n-1/2}e^{-1/4\delta^{\eta}}w_{n-\beta/2}\left(x\right) \label{eq:Remainder1}\\
 &   +D\!\!\!\!\!\!\sum_{\frac{n+1}{2}\leq i\leq n+1}\!\!\!\!\!\!\!\!\left(x+\delta^{\left(1-\eta\right)/2}\right)^{2i-n-1}\!\!\delta^{-i-1/2}\nu_{s}(0,\delta^{\left(1-\eta\right)/2})w_{n-\beta/2}(x)e^{-\frac{x^{2}}{2\delta}}\!.\nonumber 
\end{align}
The first term on the right-hand side of (\ref{eq:Remainder1}) vanishes
as $\delta\rightarrow0$, hence we have 
\begin{align}
\liminf_{\delta\rightarrow0}\mathbb{E}\int_{0}^{t}\bigl\Vert w_{n-\beta/2}R_{s,\delta}^{\left(n+1\right)}\bigr\Vert _{2}^{2}ds 
& \leq \liminf_{\delta\rightarrow0}D'\!\!\!\!\!\sum_{\frac{n+1}{2}\leq i\leq n+1}\int_{0}^{\infty}\!\!\left(x+1\right)^{4i-2n-2}x^{2n-\beta}e^{-x^{2}}dx \nonumber \\
 & \qquad \times \delta^{\left(1-\eta\right)\left(2i-n-1\right)}\delta^{-2i-1}\delta^{\frac{1}{2}\left(1-\eta\right)\left(3+\beta+\varepsilon\right)}\delta^{1/2}\delta^{n-\beta/2}\label{eq:Remainder2}
\end{align}
where we have used Proposition \ref{TDCE2} with $\beta$ replaced
by $\beta+\varepsilon<\pi/\alpha-1$ for small $\varepsilon>0$, and
$D'>0$ is a numerical constant. (Note as we integrate over $s\in\left(0,t\right)$,
the singular factor of $s^{-\gamma}$ from Proposition \ref{TDCE2}
integrates to a finite value because $\gamma\in\left(0,1\right)$.)
The $\delta$ term in (\ref{eq:Remainder2}) is bounded by 
\[
\delta^{\frac{1}{2}\left(1-\eta\right)\varepsilon-\eta\left(2n+5/2+\beta/2\right)}
\]
for small $\delta$, and, since we choose $\eta$ small enough so
that 
\[
\frac{1}{2}\left(1-\eta\right)\varepsilon>\eta\left(2n+5/2+\beta/2\right),
\]
the limit on the right-hand side of (\ref{eq:Remainder2}) vanishes.
\end{proof}


\subsection*{The inductive proof for the regularity result}

\subsubsection*{Initial case:}

The initial case is when $n=-1$, and the result we shall demonstrate
is: 
\begin{equation}
\liminf_{\delta\rightarrow0}\mathbb{E}\int_{0}^{T}\bigl\Vert w_{-2-\beta/2}T_{t,\delta}^{\left(-1\right)}\bigr\Vert _{2}^{2}dt<\infty,\quad\textrm{for every }\beta<\pi/\alpha-1.\label{eq:InitialCase}
\end{equation}
By truncating the range of integration, we have
\[
T_{t,\delta}^{\left(-1\right)}\left(x\right)  =  \int_{0}^{\infty}\int_{0}^{x}G_{\delta}\left(z,y\right)dzV_{t}\left(y\right)dy=\nu_{t}\left(0,x\right)+I_{t,\delta}\left(x\right),
\]
where $I_{t,\delta}\left(x\right):=\int_{2x}^{\infty}\int_{0}^{x}G_{\delta}\left(z,y\right)dzV_{t}\left(y\right)dy$.
From Proposition \ref{TDCE2}, it is clear that
\[
\mathbb{E}\int_{0}^{T}\left\Vert w_{-2-\beta/2}\nu_{t}\left(0,\cdot\right)\right\Vert _{2}^{2}dt<\infty
\]
and therefore we are done provided $\liminf_{\delta\rightarrow0}\mathbb{E}\int_{0}^{T}\left\Vert w_{-2-\beta/2}I_{t,\delta}\right\Vert _{2}^{2}dt<\infty$.
Using the inequality 
\[
G_{\delta}\left(z,y\right)\leq\frac{2zy}{\sqrt{2\pi\delta^{3}}}\exp\left\{ -\left(z-y\right)^{2}/2\delta\right\} ,
\]
we have a numerical constant $C>0$ such that
\[
I_{t,\delta}\left(x\right)  \leq C\delta^{-3/2}x^{2}e^{-x^{2}/2\delta}\int_{0}^{\infty}\left(y+2x\right)e^{-y^{2}/2\delta}V_{t}\left(y+2x\right)dy,
\]
and so 
\begin{multline*}
 \left\Vert w_{-2-\beta/2}I_{t,\delta}\right\Vert _{2}^{2}\\
 \leq\frac{C^{2}}{\delta^{3}}\!\!\int_{0}^{\infty}\!\!\!\!\int_{0}^{\infty}\!\!\!\!\int_{0}^{\infty}\!\!\!x^{-\beta}(y_{1}+2x)(y_{2}+2x)e^{-\frac{x^{2}+y_{1}^{2}+y_{2}^{2}}{2\delta}}V_{t}(y_{1}+2x)V_{t}(y_{2}+2x)dy_{1}dy_{2}dx.
\end{multline*}

Now let $\eta\in\left(0,1\right)$. Splitting the region of integration
on $\left(0,\delta^{\left(1-\eta\right)/2}\right)^{3}\subseteq\left(0,\infty\right)^{3}$
and its complement, we see that on the complement the triple integral
on the right-hand side above is bounded by a multiple of $e^{-1/\delta^{\eta}}$,
and so vanishes as $\delta\rightarrow0$. Therefore we need only consider
the region of integration $\left(0,\delta^{\left(1-\eta\right)/2}\right)^{3}$,
and so we have a numerical constant $C'>0$ such that 
\begin{align*}
\liminf_{\delta\rightarrow0}\mathbb{E}\int_{0}^{T}\left\Vert w_{-2-\beta/2}I_{t,\delta}\right\Vert _{2}^{2}dt
   & \leq  C'\liminf_{\delta\rightarrow0}\delta^{-\left(3+\beta\right)/2}\int_{0}^{T}\mathbb{E}\nu_{t}\left(0,\delta^{\left(1-\eta\right)/2}\right)^{2}dt\\
   & \leq  B\liminf_{\delta\rightarrow0}\delta^{\left(1-\eta\right)\left(3+\beta+\varepsilon\right)/2-\left(3+\beta\right)/2},
\end{align*}
where the second line follows by Proposition \ref{TDCE2} and $\varepsilon>0$
is sufficiently small. Clearly we may choose $\eta>0$ small enough
so that the exponent of $\delta$ is positive, whereby the right-hand
side vanishes and we have (\ref{eq:InitialCase}).

\subsubsection*{Inductive step:}

Assume that for some $-1\leq n\leq N-1$ we have 
\begin{equation}
\liminf_{\delta\rightarrow0}\mathbb{E}\int_{0}^{T}\bigl\Vert w_{k-1-\beta/2}T_{t,\delta}^{\left(k\right)}\bigr\Vert _{2}^{2}dt<\infty, \quad\textrm{for every }\beta<\pi/\alpha-1\label{eq:WhatWeWant}
\end{equation}
for $k=-1,0,1,\ldots n$. (We demonstrated the case $n=-1$ above.)
The proof will be complete if we can show condition (\ref{eq:WhatWeWant})
holds for $k=n+1$, since then, by induction, the condition holds
for all $k=0,\ldots,N$ and we can apply Lemma \ref{liminfLemma}
to give the regularity result of Theorem \ref{Reg} for the probabilistic
solution $\nu$. 

Let us begin by multiplying throughout equation (\ref{eq:Inter})
by the square of the weighting function $w_{n-\beta/2}^{2}$ and integrating
over $\left(0,\infty\right)$. To do so we apply Lemma \ref{IBP}
and arrive at
\begin{eqnarray}
 &  & \sigma_{I}^{2}\mathbb{E}\int_{0}^{t}\bigl\Vert w_{n-\beta/2}T_{s,\delta}^{\left(n+1\right)}\bigr\Vert _{2}^{2}ds \label{eq:Important}\\
 &  & \qquad=\bigl\Vert w_{n-\beta/2}T_{0,\delta}^{\left(n\right)}\bigr\Vert _{2}^{2}-\mathbb{E}\bigl\Vert w_{n-\beta/2}T_{t,\delta}^{\left(n\right)}\bigr\Vert _{2}^{2}\nonumber \\
 &  & \qquad\quad-2\mu\mathbb{E}\int_{0}^{t}\bigl\Vert w_{n-\beta/2}T_{s,\delta}^{\left(n\right)}\bigr\Vert _{2}^{2}ds\nonumber \\
 &  & \qquad\quad+2\mu\bigl(n-\frac{\beta}{2}\bigr)\mathbb{E}\int_{0}^{t}\bigl\Vert w_{n-1/2-\beta/2}T_{s,\delta}^{\left(n\right)}\bigr\Vert _{2}^{2}ds\nonumber \\
 &  & \qquad\quad+2\bigl(\sigma_{M}^{2}+\sigma_{I}^{2}\bigr)\mathbb{E}\int_{0}^{t}\bigl\Vert w_{n-\beta/2}T_{s,\delta}^{\left(n\right)}\bigr\Vert _{2}^{2}ds\nonumber \\
 &  & \qquad\quad-4\bigl(\sigma_{M}^{2}+\sigma_{I}^{2}\bigr)\bigl(n-\frac{\beta}{2}\bigr)\mathbb{E}\int_{0}^{t}\bigl\Vert w_{n-1/2-\beta/2}T_{s,\delta}^{\left(n\right)}\bigr\Vert _{2}^{2}ds\nonumber \\
 &  & \qquad\quad+2\bigl(\sigma_{M}^{2}+\sigma_{I}^{2}\bigr)\bigl(n-\frac{\beta}{2}\bigr)\bigl(n-\frac{\beta+1}{2}\bigr)\mathbb{E}\int_{0}^{t}\bigl\Vert w_{n-1-\beta/2}T_{s,\delta}^{\left(n\right)}\bigr\Vert _{2}^{2}ds\nonumber \\
 &  & \qquad\quad+\sigma_{M}^{2}\mathbb{E}\int_{0}^{t}\bigl\Vert w_{n-\beta/2}R_{s,\delta}^{\left(n+1\right)}\bigr\Vert _{2}^{2}ds \nonumber \\
 &  & \qquad\quad+2\mu\mathbb{E}\int_{0}^{t}\int_{0}^{\infty}w_{n-\beta/2}^{2}T_{s,\delta}^{\left(n\right)}R_{s,\delta}^{\left(n+1\right)}ds \nonumber 
\end{eqnarray}
Eliminating the negative terms from the right-hand side and applying
the Cauchy--Schwarz inequality gives 
\begin{eqnarray}
 &  & \sigma_{I}^{2}\mathbb{E}\int_{0}^{t}\bigl\Vert w_{n-\beta/2}T_{s,\delta}^{\left(n+1\right)}\bigr\Vert _{2}^{2}ds\label{eq:FinalReg}\\ 
 &  & \qquad\leq\bigl\Vert w_{n-\beta/2}T_{0,\delta}^{\left(n\right)}\bigr\Vert _{2}^{2} \nonumber \\
 &  & \qquad\quad-2\mu\mathbb{E}\int_{0}^{t}\bigl\Vert w_{n-\beta/2}T_{s,\delta}^{\left(n\right)}\bigr\Vert _{2}^{2}ds \nonumber\\
 &  & \qquad\quad+2\mu\bigl(n-\frac{\beta}{2}\bigr)\mathbb{E}\int_{0}^{t}\bigl\Vert w_{n-1/2-\beta/2}T_{s,\delta}^{\left(n\right)}\bigr\Vert _{2}^{2}ds\nonumber \\
 &  & \qquad\quad+2\bigl(\sigma_{M}^{2}+\sigma_{I}^{2}\bigr)\mathbb{E}\int_{0}^{t}\bigl\Vert w_{n-\beta/2}T_{s,\delta}^{\left(n\right)}\bigr\Vert _{2}^{2}ds\nonumber \\
 &  & \qquad\quad+2\bigl(\sigma_{M}^{2}+\sigma_{I}^{2}\bigr)\bigl(n-\frac{\beta}{2}\bigr)\bigl(n-\frac{\beta+1}{2}\bigr)\mathbb{E}\int_{0}^{t}\bigl\Vert w_{n-1-\beta/2}T_{s,\delta}^{\left(n\right)}\bigr\Vert _{2}^{2}ds\nonumber \\
 &  & \qquad\quad+\sigma_{M}^{2}\mathbb{E}\int_{0}^{t}\bigl\Vert w_{n-\beta/2}R_{s,\delta}^{\left(n+1\right)}\bigr\Vert _{2}^{2}ds \nonumber \\
 &  & \qquad\quad+2\left|\mu\right|\left\{ \mathbb{E}\int_{0}^{t}\bigl\Vert w_{n-\beta/2}T_{s,\delta}^{\left(n\right)}\bigr\Vert _{2}^{2}ds\right\} ^{1/2}\left\{ \mathbb{E}\int_{0}^{t}\bigl\Vert w_{n-\beta/2}R_{s,\delta}^{\left(n+1\right)}\bigr\Vert _{2}^{2}ds\right\} ^{1/2}.\nonumber 
\end{eqnarray}
We are done provided the seven terms on the right-hand side of inequality
(\ref{eq:FinalReg}) have a finite limit-infimum as $\delta\downarrow0$.
This is true of the second to fifth terms by the inductive hypothesis
(\ref{eq:WhatWeWant}), and, with this, the final two terms vanish
by Lemma \ref{RemainderTerm}. From Lemma \ref{InitialDensity}, the
first term remains finite in the limit-infimum, hence we have shown
the probabilistic solution $\nu$ satisfies the regularity properties
in Theorem \ref{Reg}. 


\subsection*{Uniqueness}

It remains to show that the probabilistic solution is indeed the unique
finite-measure valued solution of (\ref{eq:TheWeakProblem}). To do
so, suppose we have another finite-measure valued solution $\bar{\nu}$,
and define the signed-measure valued process, $\Delta$, as the difference
between the solutions: 
\[
\Delta_{t}\left(S\right):=\nu_{t}\left(S\right)-\bar{\nu}_{t}\left(S\right).
\]

Let us write $I\left(\zeta\right)$ for the expected value of a random
measure, $\zeta$, that is 
\[
I\left(\zeta\right)\left(S\right):=\mathbb{E}\left[\zeta\!\left(S\right)\right],
\]
for $S\subseteq\left(0,\infty\right)$ measurable, then, by taking
expectation in (\ref{eq:TheWeakProblem}), we have that $I\left(\nu\right)$
and $I\left(\bar{\nu}\right)$ solve 
\begin{equation}
d\left\langle \phi,\zeta_{t}\right\rangle =\frac{\sigma_{M}^{2}+\sigma_{I}^{2}}{2}\left\langle \phi'',\zeta_{t}\right\rangle dt,\quad\zeta_{0}=V_{0}.\label{eq:1DHeatEqn}
\end{equation}
Also, if we define the expected product measure to be the unique extension
of
\[
J\left(\zeta\right)\left(S_{1}\times S_{2}\right):=\mathbb{E}\left[\zeta\left(S_{1}\right)\zeta\left(S_{2}\right)\right]
\]
on $\left(0,\infty\right)^{2}$, applying It\=o's formula to the product
$\left\langle \phi_{1},\nu_{t}\right\rangle \left\langle \phi_{2},\nu_{t}\right\rangle $
using (\ref{eq:TheWeakProblem}) gives that $J\left(\nu\right)$ and
$J\left(\bar{\nu}\right)$ solve 
\begin{equation}
d\left\langle \phi,\zeta_{t}\right\rangle =\frac{\sigma_{M}^{2}+\sigma_{I}^{2}}{2}\left\langle \left\{ \partial_{xx}+2\rho^{2}\partial_{xy}+\partial_{yy}\right\} \phi,\zeta_{t}\right\rangle dt,\quad\zeta_{0}\left(x,y\right)=V_{0}\left(x\right)V_{0}\left(y\right)\label{eq:2DHeatEqn}
\end{equation}
for $\phi\in C_{0}^{\infty}\left(\left(0,\infty\right)^{2}\right)$.
Equations (\ref{eq:1DHeatEqn}) and (\ref{eq:2DHeatEqn}) are just
the deterministic heat equations with zero boundary conditions, hence
their solutions are unique so we conclude $I\left(\nu\right)=I\left(\bar{\nu}\right)$
and $J\left(\nu\right)=J\left(\bar{\nu}\right)$ --- that is, the
first and second moments of $\nu$ and $\bar{\nu}$ are equal in distribution.

By the linearity of (\ref{eq:TheWeakProblem}), we have that the smoothed
measure \[T_{\delta}\Delta_{t}=\int_{0}^{\infty}G_{\delta}\left(\cdot,y\right)\Delta_{t}\left(dy\right)\]
solves 
\[
dT_{\delta}\Delta_{t}=-\mu\left(T_{\delta}^{r}\Delta_{t}\right)'dt+\frac{1}{2}\left(\sigma_{M}^{2}+\sigma_{I}^{2}\right)\left(T_{\delta}\Delta_{t}\right)''dt-\sigma_{M}\left(T_{\delta}^{r}\Delta_{t}\right)'dW_{t}
\]
with $T_{t,\delta}\Delta\equiv0$. Using the work in (\ref{eq:Inter})
and integrating over $\Omega\times\left(0,\infty\right)$ gives 
\[
\mathbb{E}\left\Vert T_{\delta}\Delta_{t}\right\Vert _{2}^{2}
  \leq  \int_{0}^{t}\sigma_{M}^{2}\mathbb{E}\left\Vert \left(R_{s,\delta}\Delta\right)'\right\Vert _{2}^{2}ds+2\left|\mu\right|\int_{0}^{t}\mathbb{E}\bigl[\left\Vert T_{\delta}\Delta_{s}\right\Vert _{2}^{2}\bigr]^{\frac{1}{2}}\mathbb{E}\bigl[\left\Vert \left(R_{s,\delta}\Delta\right)'\right\Vert _{2}^{2}\bigr]^{\frac{1}{2}}ds,
\]
then, since 
\[
\mathbb{E}\left[\left|\Delta_{t}\left(0,\varepsilon\right)\right|^{2}\right]\leq\mathbb{E}\left[\left|\nu_{t}\left(0,\varepsilon\right)\right|^{2}\right]+\mathbb{E}\left[\left|\bar{\nu}_{t}\left(0,\varepsilon\right)\right|^{2}\right]=2\mathbb{E}\left[\nu_{t}\left(0,\varepsilon\right)^{2}\right],
\]
Lemma \ref{RemainderTerm} can be applied to $R_{s,\delta}\Delta$
and gives that $\liminf_{\delta\rightarrow0}\mathbb{E}\left[\left\Vert T_{\delta}\Delta_{t}\right\Vert _{2}^{2}\right]=0$.
We conclude that we have uniqueness, since, for every $\phi\in L^{2}\left(0,\infty\right)$,
\begin{align*}
\mathbb{E}\int_{0}^{T}\left(\left\langle \phi,\nu_{t}\right\rangle -\left\langle \phi,\bar{\nu}_{t}\right\rangle \right)^{2}dt 
  & =  \mathbb{E}\int_{0}^{T}\lim_{\delta\rightarrow0}\left\langle \phi,T_{\delta}\Delta_{t}\right\rangle ^{2}dt \\
  & \leq \left\Vert \phi\right\Vert _{2}^{2}\int_{0}^{T}\liminf_{\delta\rightarrow0}\mathbb{E}\left\Vert T_{\delta}\Delta_{t}\right\Vert _{2}^{2}dt=0,
\end{align*}
and so, by approximation, the result holds for every $\phi\in C^{\mathrm{test}}$.
\qed


\section{Proof of Theorem \ref{Converse}}

This theorem relies on the fact that the estimate in Proposition \ref{TDCE2}
is sharp. Specifically, Lemma 5.1 of \cite{Jin} states that if $V_{0}$
is supported on a finite interval away from zero, then for every $\eta>0$
there exist $C\left(\eta\right)>0$ and $\varepsilon_{0}\left(\eta\right)>0$
such that
\begin{equation}
\mathbb{E}\left[\nu_{t}\left(0,\varepsilon\right)^{2}\right]\geq C\left(\eta\right)\varepsilon^{2+\pi/\alpha+\eta},\quad\mathrm{for\, all\,}0<\varepsilon<\varepsilon_{0}\left(\eta\right)\,\mathrm{and}\, t\in\left[0,T\right].\label{eq:JinIne}
\end{equation}
It is clear that if $\tilde{V}_{0}\leq V_{0}$, then the empirical
measure, $\tilde{\nu}$, obtained from starting the system from $\tilde{V}_{0}$,
satisfies $\tilde{\nu}_{t}\leq\nu_{t}$. Since we can certainly find
a constant $a>0$ such that $\tilde{V}_{0}:=V_{0}\mathbf{1}_{\left(a,\infty\right)}$
is non-zero, we see that (\ref{eq:JinIne}) remains valid if we drop
the assumption on the support of $V_{0}$. 

To complete the proof, suppose, for a contradiction, that we have
$\beta>\pi/\alpha-1$ and a $k\in\left\{ 0,1,\ldots,N\right\} $ such
that
\[
\mathbb{E}\int_{0}^{T}\bigl\Vert w_{k-1-\beta/2}V_{t}^{\left(k\right)}\bigr\Vert _{2}^{2}dt<\infty.
\]
Then 
\[
\mathbb{E}\int_{0}^{T}\int_{0}^{1}x^{2\left(k-1\right)-\tilde{\beta}}V_{t}^{\left(k\right)}\left(x\right)^{2}dxdt<\infty
\]
for any $\tilde{\beta}<\beta$, and let us fix a $\tilde{\beta}\in\left(\pi/\alpha-1,1\right)$
and take $\varepsilon>0$ small enough so that we have $\tilde{\beta}-\varepsilon\in\left(\pi/\alpha-1,1\right)$.
Then for $k\geq2$ 
\[
\int_{0}^{1}x^{2\left(k-2\right)-\tilde{\beta}+\varepsilon}dx<\infty,
\]
and, since 
\begin{align}
\Bigl(V_{t}^{\left(k-1\right)}\left(x\right)-V_{t}^{\left(k-1\right)}\left(1\right)\Bigr)^{2} 
  & = \Bigl(\int_{x}^{1}V_{t}^{\left(k\right)}\left(y\right)dy\Bigr)^{2}\label{eq:SwitchMe}\\
  & \leq \int_{0}^{1}x^{2\left(k-1\right)-\tilde{\beta}}V_{t}^{\left(k\right)}\left(x\right)^{2}dx\cdot\int_{x}^{1}y^{-2k+2+\tilde{\beta}}dy\nonumber \\
  & \leq \int_{0}^{1}x^{2\left(k-1\right)-\tilde{\beta}}V_{t}^{\left(k\right)}\left(x\right)^{2}dx\cdot\frac{x^{-2k+3+\tilde{\beta}}}{2k-3-\tilde{\beta}},\nonumber 
\end{align}
applying the triangle inequality, multiplying by $x^{2k-4-\tilde{\beta}+\varepsilon}$,
and integrating over $\Omega\times\left[0,T\right]\times\left(0,1\right)$
gives
\[
\mathbb{E}\int_{0}^{T}\int_{0}^{1}x^{2\left(k-2\right)-\tilde{\beta}+\varepsilon}V_{t}^{\left(k-1\right)}\left(x\right)^{2}dxdt<\infty.
\]
Repeating this argument eventually yields 
\[
\mathbb{E}\int_{0}^{T}\int_{0}^{1}x^{-\tilde{\beta}+\varepsilon}V_{t}'\left(x\right)^{2}dxdt<\infty
\]
for all allowable $\varepsilon$. If we return to (\ref{eq:SwitchMe})
with $k=1$ and take the upper limit of the integral to be $x$ and
the lower limit $0$, then, since $V_{t}\left(0\right)=0$, we arrive
at
\[
\mathbb{E}\int_{0}^{T}\int_{0}^{1}x^{-2-\tilde{\beta}+\varepsilon}V_{t}\left(x\right)^{2}dxdt<\infty.
\]
Repeating this manipulation once more gives the following for $x<1$:
\begin{align*}
\mathbb{E}\int_{0}^{T}\nu_{t}\left(0,x\right)^{2}dt & \leq  \mathbb{E}\int_{0}^{T}\int_{0}^{1}x^{-2-\tilde{\beta}+\varepsilon}V_{t}\left(x\right)^{2}dxdt\int_{0}^{x}y^{2+\tilde{\beta}-\varepsilon}dy\\
 & =  D(\tilde{\beta},\varepsilon)x^{3+\tilde{\beta}-\varepsilon},
\end{align*}
with $D(\tilde{\beta},\varepsilon)\in\left(0,\infty\right)$
a constant. This is a contradiction of (\ref{eq:JinIne}) for the
case $\eta=(1-\pi/\alpha)+(\tilde{\beta}-\varepsilon)>0$,
and hence we have Theorem \ref{Converse}. \qed


{\setstretch{1.0}
  \bibliographystyle{alpha}
  \bibliography{reference}

\newcommand{\etalchar}[1]{$^{#1}$}
\begin{thebibliography}{BHH{\etalchar{+}}11}

\bibitem[BC09]{Bain}
A.~Bain and D.~Crisan.
\newblock {\em Fundamentals of Stochastic Filtering}.
\newblock Springer, New York, 2009.

\bibitem[BHH{\etalchar{+}}11]{Bush}
N.~Bush, B.M. Hambly, H.~Haworth, L.~Jin, and C.~Reisinger.
\newblock {Stochastic Evolution Equations in Portfolio Credit Modelling}.
\newblock {\em SIAM Journal of Financial Mathematics}, 2(1):627--664, 2011.

\bibitem[BL95]{Bichteler}
K.~Bichteler and S.J. Lin.
\newblock {On the stochastic Fubini theorem}.
\newblock {\em Stochastics and Stochastic Reports}, 54(3-4):271--279, 1995.

\bibitem[Iye85]{Iyengar}
Satish Iyengar.
\newblock Hitting lines with two-dimensional brownian motion.
\newblock {\em SIAM Journal of Applied Mathematics}, 45(6):983--989, 1985.

\bibitem[Jin10]{Jin}
L.~Jin.
\newblock {\em Particle systems and {SPDE}s with applications to credit
  modelling}.
\newblock D.Phil {T}hesis, {U}niversity of {O}xford, 2010.

\bibitem[KK04]{Kim}
K-H. Kim and N.~V. Krylov.
\newblock {On SPDEs with Variable Coefficients in One Space Dimension}.
\newblock {\em Potential Analysis}, 21(3):209--239, 2004.

\bibitem[KL98]{KryLoto}
N.V. Krylov and S.V. Lototsky.
\newblock {\em {A Sobolev space theory of SPDEs with constant coefficients in a
  half line}}, volume~30.
\newblock 1998.

\bibitem[Kry94]{Krylov}
N.V. Krylov.
\newblock A ${W}_{2}^{n}$-theory of the {D}irichlet problem for {SPDE}s in
  general smooth domains.
\newblock {\em Probability Theory and Related Fields}, 98(3):389--421, 1994.

\bibitem[Kry03]{KrylovHeat}
N.V. Krylov.
\newblock Brownian trajectory is a regular lateral boundary for the heat
  equation.
\newblock {\em {SIAM J. Math. Anal.}}, 34(5):1167--1182, 2003.

\bibitem[KX99]{Kurtz}
T.G. Kurtz and J.~Xiong.
\newblock Particle representations for a class of nonlinear {SPDEs}.
\newblock {\em Stochastic Processes and their Applications}, 83(1):103--126,
  1999.

\bibitem[Met10]{Metzler}
A.~Metzler.
\newblock On the first passage problem for correlated {B}rownian motion.
\newblock {\em Statistics \& Probability Letters}, 80(5-6):277--284, 2010.

\end{thebibliography}
	
}
\end{document}